\documentclass[11pt]{article}

\usepackage[english]{babel}
\usepackage[utf8]{inputenc}
\usepackage[T1]{fontenc}
\usepackage{lmodern}
\usepackage[a4paper, margin=2.2cm]{geometry}
\usepackage{graphicx}
\usepackage{booktabs}
\usepackage{caption}
\usepackage{microtype}
\usepackage[section]{placeins}
\usepackage{mathtools}
\usepackage{amssymb}
\usepackage{csquotes}
\usepackage{tikz}
\usetikzlibrary{calc} % REQUIRED by the KAYROS wordmark coordinates below

\PassOptionsToPackage{hyphens}{url}
\usepackage[
  backend=biber,
  style=alphabetic,
  sorting=none,
  language=american,
  maxcitenames=1,
  mincitenames=1,
  maxbibnames=7,
  minbibnames=7,
  alldates=ymd
]{biblatex}
\usepackage[
  colorlinks=true,
  linkcolor=blue,
  citecolor=blue,
  urlcolor=blue
]{hyperref}
\hypersetup{
  pdftitle={Thread Scaling of Hexaly on the TDVRPTW across Two Model Encodings. An Experimental Report: Binding Choice, Slice-Count Choice, and Two Thread Ladders},
  pdfauthor={Florian Rascoussier}
}
\usepackage[capitalise, noabbrev]{cleveref}

\newcommand{\keywords}[1]{\textbf{Keywords}: \textit{#1}}
\newcommand{\sgap}{\ensuremath{\varrho}}
\newcommand{\ascore}{\ensuremath{\Lambda}}

\newcommand{\SlashedO}{%
  \tikz[baseline=(char.base)]{
    \node[inner sep=0pt] (char) {O};
    \begin{scope}
    \clip (char.south west) rectangle (char.north east);
    \draw[line width=.55pt]
         ($(char.north west)+(2.1pt,-1pt)$) -- ($(char.south east)+(-2.1pt,1pt)$);
    \end{scope}
}}

\DeclareRobustCommand{\kayrosinner}{KAYR\SlashedO S}

\newcommand{\kayros}{\texorpdfstring{\kayrosinner}{KAYROS}}

\title{\Large Thread Scaling of Hexaly on the TDVRPTW across Two Model Encodings\\[1.2ex]
\large An Experimental Report: Binding Choice, Slice-Count Choice, and Two Thread Ladders}
\author{Florian Rascoussier\thanks{IMT Atlantique, Lab-STICC, CNRS, UMR 6285 (équipe DECIDE) and INSA Lyon, Inria, CITI, UR3720, 69621 Villeurbanne, France (équipe EMERAUDE); ORCID identifier: \href{https://orcid.org/0009-0005-3253-9814}{0009-0005-3253-9814}\@.}}
\date{Version 3, 19/10/2026}

\begin{document}

\maketitle

\begin{abstract}
Thread scaling in Hexaly on the Time-Dependent Vehicle Routing Problem with Time Windows (TDVRPTW) depends on modeling, and on the language binding through which the model reaches the solver. We compare two encodings, both reaching the solver through its C++ binding: one evaluates continuous travel-time functions exactly through external callbacks, which that binding evaluates concurrently, while the other approximates them with time slices evaluated natively by the solver. Independent CPU accounting confirms that both use the cores they are allocated, but they respond differently to width. The external-function encoding stabilizes, reducing its seed dispersion by about 30\% up to 8 threads before widening again, yet its pooled quality barely moves, because nine of the ten instances improve while the hardest one degrades by enough to cancel them. For the time-sliced encoding we first select a discretization on feasibility under the original travel-time functions rather than on approximation error, and the retained setting then improves monotonically with the thread count, gaining about a quarter of its final gap and reducing the seed dispersion by nearly half between 1 and 16 threads. At equal single-threaded budget the two encodings are close, and the small pooled edge of exact evaluation comes from that same hardest instance, so the practical advantage of the discretization is thread scaling rather than fidelity traded for speed. These descriptive results support multi-threading as a way to improve solution quality and stability when the model can exploit it, and they show that binding-level constraints on external evaluation must be measured, not assumed.

\vspace{1ex}
\end{abstract}

\begin{center}
\begin{minipage}{0.92\linewidth}
\small \textbf{Note on this version.} Version 3 updates the author lists in the references to KAYROS, its technical report, MAMUT-routing and Poryos2026. It adds no experiments and preserves the results and conclusions of version 2. The first version of this report measured the external-function encoding through the Python binding of the solver, which declares external functions as callbacks from a Global Interpreter Lock language and therefore evaluates them on a single worker at every requested width. That version consequently described the encoding as unable to scale, and reported the time-sliced encoding as uniformly superior. Version 2 replaced the entire external-function ladder with runs of the same model through the C++ binding, allowing concurrent evaluation, and ported the time-sliced arm to the same discipline: its thread ladder and the slice-count study were rerun through the C++ binding at identical cells and seeds. For the external-function encoding the two bindings are shown to evaluate identical values. For the time-sliced one a per-run parity check shows that both bindings build bit-identical sliced tables. Both substitutions therefore repair a measurement of the interface rather than changing the model. The results of both encodings, their comparison at equal budget, and the conclusions that rested on them were revised accordingly. The first version remains the record of what the Python binding delivers.
\end{minipage}
\end{center}

\clearpage
\setcounter{tocdepth}{2}
\tableofcontents
\clearpage

\section{Context and scope}
\label{sec:context}

This report belongs to a solver-comparison campaign on the duration-minimization TDVRPTW, carried out as part of PhD research that also produces an open-source solver, \kayros{} \autocite{rascoussierKAYROS2026,rascoussierKAYROSTechReport2026}. Hexaly Optimizer \autocite{hexalyOptimizer2026} is one of the solvers considered in this campaign, under an academic license that is gratefully acknowledged. The main comparison assigns the same wall-clock budget to every solver and runs each process single-threaded on one physical core. This convention equalizes resources rather than constraining any particular solver: every arm of the campaign, \kayros{} included, receives one core and one hour per cell.

During the exchange that motivated this report, Hexaly explained that their optimizer is designed for multi-threaded use and that restricting it to one thread may reduce performance and solution quality. That objection called for a dedicated thread-scaling study, conducted here under a similar instance subset and same wall-clock budget as the main comparison. Maxime Rougier also clarified that Hexaly primarily uses additional threads to diversify the search and the methods applied. Stability across seeds is therefore examined alongside mean solution quality.

The outcome depends on the model encoding and on the binding through which the model reaches the solver, and the report is organized around this distinction. \Cref{sec:protocol} introduces the experimental protocol and the two encodings. \Cref{sec:tier1} measures the thread scaling of the external-function encoding through the C++ binding, and explains why the Python binding used in the first version of this report could not measure it. \Cref{sec:slices} selects a sufficiently fine time discretization, \cref{sec:tier2} reports the thread scaling of the time-sliced encoding and compares the two encodings at equal budget, and \cref{sec:conclusion} draws the implications.

\section{Experimental protocol}
\label{sec:protocol}

\subsection{Instances and reference values}
\label{sec:instances}

The thread ladders use 10 TDVRPTW instances from five benchmark families, spanning small literature-derived cases to large urban-traffic instances (\cref{tab:instances}). The panel draws from Dabia2013 \autocite{dabiaBranchPriceTimeDependent2013,solomonAlgorithmsVehicleRouting1987,ichouaVehicleDispatchingTimedependent2003}, Rifki2020 \autocite{rifkiImpactSpatiotemporalGranularity2020}, Vu2020 \autocite{vuDynamicDiscretizationDiscovery2020}, our Lera2026 family derived from the Gehring--Homberger bases \autocite{hombergerTwoEvolutionaryMetaheuristics1999}\footnote{Lera2026 is named in honor of Gonzalo Lera-Romero, whose open-source exact TDVRPTW solver \autocite{lera-romeroLinearEdgeCosts2020} shaped our work. Lera-Romero is not an author of the family.}, and our OpenStreetMap-based Poryos2026 family \autocite{rascoussierPoryos2026}. The instances and their piecewise-linear travel-time data are distributed as checksummed artifacts through MAMUT-routing \autocite{pichon:hal-05629810v1}.

\begin{table}[ht]
\centering
\caption{The ten-instance set of the thread ladders. The handle column gives the short name used in the result tables. Reference values are checker-validated best-known solutions, or checker-exact optimality certificates where marked. The slice study of \cref{sec:slices} uses a six-instance subset: the five marked with a dagger, plus RC207 (Dabia2013, $n = 100$).}
\label{tab:instances}
\small
\begin{tabular}{llllr}
\toprule
Handle & Family & Instance & Reference & $n$ \\
\midrule
R109 & Dabia2013 & R109 & certificate & 50 \\
RC103$^{\dagger}$ & Dabia2013 & RC103 & certificate & 100 \\
Rifki-27$^{\dagger}$ & Rifki2020 & Rifki-27 & best known & 50 \\
Rifki-12 & Rifki2020 & Rifki-12 & best known & 60 \\
Vu-A5 & Vu2020 & Vu-A5-pB-d70-w120 & certificate & 59 \\
Vu-A4$^{\dagger}$ & Vu2020 & Vu-A4-pB-d90-w150 & certificate & 99 \\
Lera-C2 & Lera2026 & Lera-C2\_2\_5-S3 & best known & 200 \\
Lera-RC1$^{\dagger}$ & Lera2026 & Lera-RC1\_10\_5-S1 & best known & 1000 \\
SF & Poryos2026 & poryos-san\_francisco-n100-hyb-bpr-light & best known & 100 \\
Tokyo$^{\dagger}$ & Poryos2026 & poryos-tokyo-n1000-hyb-bpr-heavy & best known & 1000 \\
\bottomrule
\end{tabular}
\end{table}

The panel was frozen after exploratory observations on a larger, unreported set. It was deliberately selected to remain computationally feasible while preserving diversity in family, size, reference status, and origin. We use the same panel for both encodings to preserve cell-by-cell pairing. The results are therefore descriptive of this panel rather than estimates of an average effect over TDVRPTW instances. All gaps are measured against a reference snapshot frozen before the campaign. Since this reference set evolves with our experiments, the absolute gaps are working measurements rather than final publication values, but comparisons across thread levels remain valid because every run uses the same snapshot.

\subsection{The two Hexaly model encodings}
\label{sec:model}

Both encodings formulate the TDVRPTW with list decision variables in the Hexaly modeling API, following the general structure of Hexaly's TDCVRPTW template \autocite{hexalyTDCVRPTW2026}. Route sequences satisfy a partition constraint, vehicle capacity is imposed as a hard constraint, and customer and depot time windows are represented through a first-level total-lateness objective. The solver therefore minimizes the lexicographic pair (total lateness, total duration), with a free departure time from the depot for each route. The two encodings differ in their evaluation of time-dependent travel times, and this difference determines whether the search can use several threads.

\paragraph{The external-function encoding (Tier 1).}
Arrival times are computed through external functions that evaluate the same piecewise-linear arrival-time data, with the same conventions, as our reference checker. This removes the travel-time discretization error of the second encoding, but the objectives are not equivalent because the checker independently optimizes depot departure times. Their differences are small but nonzero, and checker values are authoritative throughout this report. The callback is side-effect free: it reads an immutable table while all other variables are local. The encoding exists in two implementations, one declaring its functions through the Python API and one through the C++ API, and the ladder of \cref{sec:tier1} uses the second. An equivalence gate ties them together and demonstrated that the two implementations evaluate the same model. The choice between them is a question of interface rather than of formulation.

\paragraph{The time-sliced encoding (Tier 2).}
Let $\mathcal{V}$ be the vertex set of the instance, i.e. the $n$ customers together with a single depot vertex, so that $\left| \mathcal{V} \right| = n + 1$, and let $\alpha_{ij}$ be the original arrival-time function of arc $\langle i,j \rangle$, giving the arrival time at $j$ for a departure from $i$ at a given time. Let $[h_0,h_1]$ be the travel-time horizon, let $\psi \in \mathbb{N}_{>0}$ be the \emph{slice count}, i.e. the number of equal periods into which that horizon is divided, and let $w=(h_1-h_0)/\psi$ be the resulting slice width. For arc $\langle i,j \rangle$ and slice $k \in [0, \psi-1]_\mathbb{N}$, the model stores the original travel time at the slice midpoint $\mu_k=h_0+(k+1/2)w$, namely $\alpha_{ij}(\mu_k)-\mu_k$, in a constant $\left| \mathcal{V} \right| \times \left| \mathcal{V} \right| \times \psi$ table. The midpoint evaluation follows the checker's convention of taking the smallest ordinate at a duplicate piecewise-linear abscissa. A departure time $t_\text{dep}$ propagated along a route is rounded to the nearest integer with the ties-to-even convention, then assigned to slice $\min\{\psi-1,\max\{0,\lfloor(t_\text{dep}-h_0)/w\rfloor\}\}$. An exact internal boundary therefore enters the following slice, while $h_1$ enters the final slice. A propagated departure outside the supported integer-time horizon has no defined lookup and the corresponding move is rejected. The travel time is frozen at the departure slice: the encoding neither accounts for crossing a boundary during traversal nor modifies the table to enforce FIFO, so adjacent slices can invert arrival order. Time windows, service times, demands and capacity retain their original values.

The model follows the Hexaly's published TDCVRPTW template \autocite{hexalyTDCVRPTW2026} and uses one list decision variable for each vehicle available to the solver. When an instance declares a fleet size, that value determines the number of list variables. For Poryos2026, whose files do not declare a fleet size, we allow up to $n$ vehicles and therefore create $n$ list variables. We remove the fleet-size level of the template's lexicographic objective and replace its distance level with total route duration, bringing the native objective closer to the checker duration used in the campaign. Each route's departure remains a continuous variable in the intersection of the travel-time horizon and depot time window instead of being fixed at the beginning of the horizon. Since the slice table approximates the original travel-time functions and does not enforce the first-in-first-out (FIFO) property, the slice count $\psi$ is a substantive modeling parameter studied in \cref{sec:slices}. Like the external-function encoding, this encoding exists in a Python and a C++ implementation, the latter of which is used in this report.

We discretize travel times rather than speeds. For families generated from piecewise-constant speed profiles, integrating a discretized speed profile would preserve FIFO and can represent the original model exactly. It would, however, require evaluating a departure-dependent integral inside the search, precisely the operation that the native encoding must avoid. Freezing travel time at the departure slice retains a single native table lookup but introduces boundary-crossing errors and may violate FIFO. The slice count therefore controls model fidelity and is a modeling decision rather than a tuning detail.

Hexaly 15.0 does not expose a best-solution event. We therefore poll the incumbent trajectory through a periodic callback whose tick period has a one-second lower bound, giving timestamps with an approximate one-second granularity for both encodings. We record a new incumbent whenever the native lexicographic objective has strictly improved since the previous tick.

\subsection{Execution grid, pinning, and measurement}
\label{sec:grid}

Each thread-scaling experiment combines the 10 instances with requested thread counts $\xi \in \{1, 2, 4, 8, 16\}$ and 10 seeds, for a total of 500 one-hour runs per ladder. The seeds are fixed once and reused by every arm of the study, so that any cell of one ladder pairs with the same cell of the other. The first version of this report used only three seeds, and this version adds the remaining seven for improved statistical significance. The wall-clock budget includes data loading and model construction as well as the search itself. We stop at 16 threads to cover a meaningful scaling range while keeping the campaign computationally manageable. Four ladders were run in total, a C++ and a Python one per encoding, at identical cells. The external-function Python ladder was not extended beyond the original three seeds, since its purpose is the binding comparison of \cref{sec:tier1-setup} and its verdict does not depend on the seed grid. The time-sliced Python ladder covers the full ten-seed grid and is superseded cell for cell by its C++ counterpart, which reclaims the Python-side setup time on the largest instances without changing the model. This report presents the two C++ ladders and uses the Python runs only for the binding comparisons.

Each run starts one Hexaly process, pinned with \texttt{taskset} to its requested number of physical cores on the same 32-core node class. Core sets remain exclusive to a process, but processes can still share node-level resources such as caches and memory bandwidth. Since narrow configurations allow more processes to share a node than wide configurations, thread width remains partly confounded with node co-tenancy. We account for this limitation when interpreting the native thread ladder and use an independent replication to bound its likely effect (\cref{sec:tier2-setup}).

We independently measure each complete process tree with GNU~time, recording wall time, user and system CPU time, and peak memory use. Comparing CPU time with wall time reveals how many cores the process actually keeps busy, which is essential because the solver API reports the requested thread count rather than effective parallelism. We therefore base the resource analysis on measured process consumption rather than scheduler allocation, which also includes idle and fragmented cores.

\subsection{Scoring and metrics}
\label{sec:metrics}

Anytime solver performance is reported through two complementary measures computed from the checker-valid incumbent trajectory: the final gap at the time limit and the Normalized Signed Primal Integral over the full run. The final gap describes where the search ends, while the primal integral rewards methods that find good feasible solutions early and keep improving them throughout the search. Together, they provide a balanced view of final and anytime quality. Lower values are better for both. We also report the time to the first checker-valid incumbent, and we describe dispersion by the sample standard deviation of the final gap across the ten seeds of a cell, computed per instance with the usual sample-size-minus-one denominator and then averaged over instances with equal weight. We prefer it to the max-minus-min seed range used in the first version of this report, because a range grows with the number of seeds by construction and therefore cannot be compared across seed grids, whereas a standard deviation can.

The reported trajectory is the best-so-far envelope of sampled routes accepted by the MAMUT-routing reference checker (\texttt{mamut-\allowbreak routing-\allowbreak lib} version 0.9.0) on the original travel-time functions. It contains routes captured after strict improvements in Hexaly's native objective, plus the terminal route when it is checker-valid. Rejected routes contribute nothing. The envelope is therefore checker-valid by construction, but it remains a sample of the states visited by the solver. Adding the terminal route does not change any score in these experiments because it never improves the preceding envelope. We nevertheless report terminal validity where relevant and reserve native objective values for diagnostics.

At the deadline $\Gamma = 3600$\,s, we measure end-of-run quality with the squeezed gap $\sgap(z_{\Gamma}, z^{*}) = (z_{\Gamma} - z^{*})/(z_{\Gamma} + z^{*})$, where $z_{\Gamma}$ is the best checker-valid value and $z^{*}$ is the reference duration. Tables and figures call this quantity the \emph{final gap}. It is not necessarily the gap of Hexaly's terminal native solution. The value is zero when the run matches the reference and negative when it improves upon it. We set it to 1 when no checker-valid route was recorded. This bounded measure is a rescaling of the classic relative gap $g(z) = (z - z^{*})/z^{*}$:
\[
\sgap = \frac{g}{g + 2}, \qquad g = \frac{2\sgap}{1 - \sgap},
\]
so the two gaps contain the same information. The bounded $\sgap$ scale is suitable for aggregation across instances, while the final-gap figures also provide the familiar percentage scale $g$ as an alternative axis.

The classic gap is not used as the primary scale because three of its properties are exercised directly by this campaign. It has no finite value when a run records no checker-valid incumbent, which happens at the coarse slice counts of \cref{sec:slices}, so pooling in $g$ would require an external convention such as discarding the cell or charging it an arbitrary penalty, whereas $\sgap = 1$ is the intrinsic worst value of the bounded scale and needs no such decision. It is unbounded above, so a single hard cell can dominate a mean taken over instances of very different difficulty. Finally, it treats improvement and degradation asymmetrically: a solution twice as long as the reference gives $g = 1$ while one half as long gives $g = -0.5$, where $\sgap$ maps reciprocal ratios to opposite values. The first property also governs the anytime metric below, which integrates over the period preceding the first incumbent, where the classic gap has no value to contribute.

Anytime quality is measured by the Normalized Signed Primal Integral
\[
\ascore(\Gamma) = \frac{1}{\Gamma} \int_{0}^{\Gamma} \sgap\bigl(z(t), z^{*}\bigr) \, \mathrm{d}t,
\]
where $z(t)$ is the best checker-valid incumbent available at time $t$. Before the first valid incumbent, the gap is held at 1. Both performance measures therefore belong to $(-1, 1]$. The alternative classic-gap axis of the figures is reserved for the final-gap panels, which carry a gap of the kind $g$ measures. It is not offered for $\ascore$: although the anytime metric is itself a level on the $\sgap$ scale, labelling it in $g$ would invite reading it as a time-averaged classic gap, and no such average exists, since the gap is held at 1 before the first valid incumbent, where $g$ is unbounded. For dispersion, the nonlinear gap transformation admits no unique conversion either, since a standard deviation is not recoverable from two endpoints the way a range is. The anytime and dispersion figures therefore stay on the $\sgap$ scale alone.

\subsection{Hardware and software}
\label{sec:hardware}

All experiments run on the grvingt cluster of the Grid'5000 testbed at Nancy \autocite{grid5000Testbed2026}. We use its dual-socket 32-core Intel Xeon Gold 6130 nodes with 192\,GiB of RAM. The software stack is Hexaly Optimizer 15.0, driven through its C++ API, compiled with GCC 13, for both encodings. The Python-binding runs kept for the binding comparisons used CPython 3.13.

\section{The external-function encoding: evaluation scales once the binding permits it}
\label{sec:tier1}

\subsection{Setup}
\label{sec:tier1-setup}

This first ladder measures the exact encoding through the C++ external-function binding. It uses the instances, seeds, thread levels, wall-clock budget, hardware class, execution protocol and scoring of \cref{sec:instances,sec:grid,sec:metrics}, so every cell pairs with the corresponding cell of the time-sliced ladder of \cref{sec:tier2}. All 500 runs reach the deadline, record a checker-valid incumbent, and end without lateness.

The first version of this report ran the same model through the Python binding and observed one busy evaluator at every requested width. That serialization is a property of the binding rather than of the encoding. Python external functions are declared to the optimizer as callbacks from a Global Interpreter Lock language, and the optimizer then limits the effective search width to avoid harmful concurrent access \autocite{hexalyExternalFunctions2026}, whereas C++, C\# and Java callbacks may be evaluated concurrently once the caller guarantees thread safety \autocite{hexalyCppExternalFunction2026}. Our callback only reads an immutable table, so porting it required no change of semantics, and an internal validation gate confirms that both implementations return identical values. What changes is throughput: at $n = 100$ and one thread, the C++ binding performs around 90\,000 search iterations per second against 5\,000 through the Python one, a factor of 18 (factors between 14 and 17 on the smaller gate instances). Diagnostic attempts to force parallel evaluation from Python, under a stock as well as a free-threaded interpreter, either collapsed in throughput or crashed, and never produced a supported parallel configuration. The parameter API reports the requested width rather than the effective one \autocite{hexalyCppParam2026}, so the ladder below still rests on independent CPU accounting, and the measurements of the first version remain the record of what the Python binding delivers. % removed detail: The overhead removed is interpreter dispatch, a large share of it spent on the undefined-operand probes that make up about a fifth of the calls, each of which raised and caught an exception in Python.

\subsection{CPU occupancy per thread level}
\label{sec:tier1-cpu}

\begin{figure}[ht]
\centering
\includegraphics[width=\linewidth]{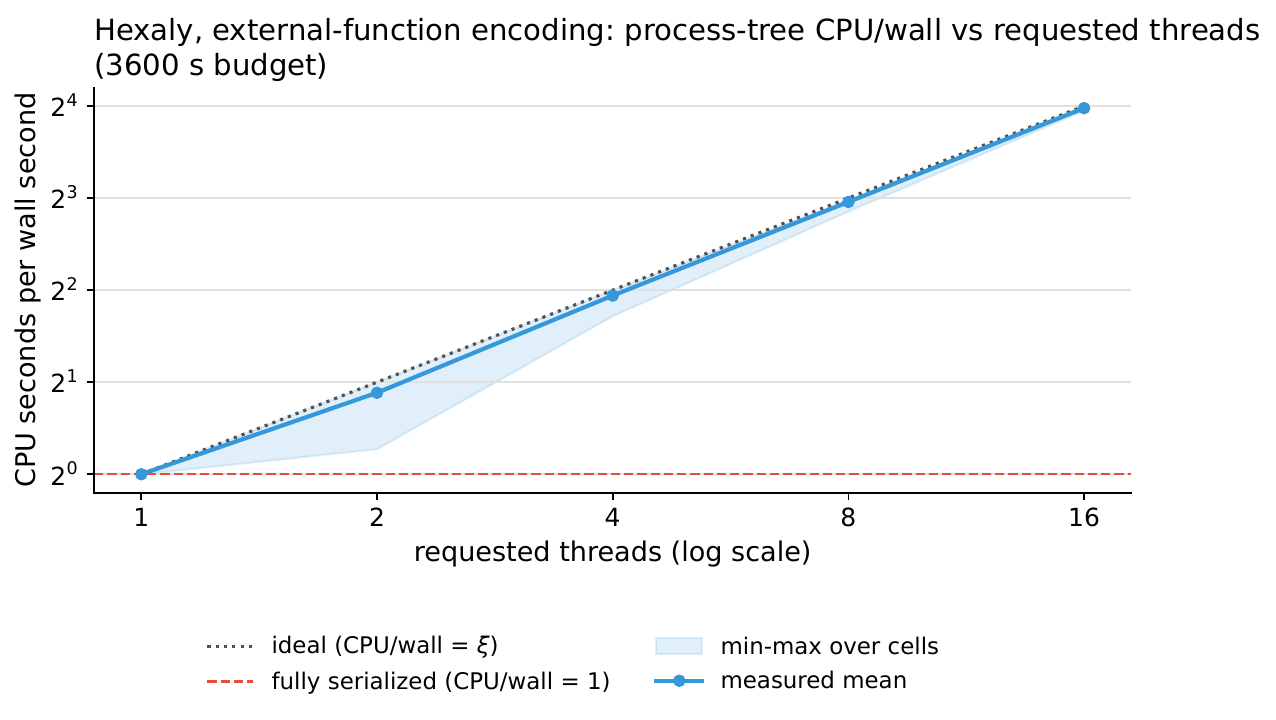}
\caption{External-function encoding under the C++ binding: measured process-tree CPU per wall-clock second against the requested thread count (GNU~time, whole process tree, 100 runs per level, mean and min--max band). The upper reference line denotes full occupancy of the allocated cores and the lower one denotes fully serialized execution, which is where the same measurement stood at every width under the Python binding. The figure shows occupancy only and says nothing about how usefully the cores are spent.}
\label{fig:tier1-cpu}
\end{figure}

The encoding now occupies the cores it is given. Measured process-tree CPU per wall-clock second follows the requested width at every level (\cref{fig:tier1-cpu}), rising from 1.00 at one thread to 15.75 at 16. Under the Python binding, the same measurement stayed at one busy core, illustrating the binding's serialization. The C++ binding therefore allows the encoding to scale, and the CPU accounting confirms that it does so.

The dip at 2 threads deserves a word, because it is the only level where the independent check flags any cell. Eleven of the 100 runs at that width report between 1.21 and 1.50 CPU seconds per wall second instead of the 1.90 that the other 89 average, and no run at any other width is flagged. All eleven come from a single batch on two nodes of the cluster, and the runs concerned are not systematically worse than their unflagged siblings. A dedicated re-run of the eleven cells, executed alone on one host with no other job sharing it, recovers full occupancy in every one of them, each measuring between 1.999 and 2.000 CPU seconds per wall second, which places the effect in the packing of that batch rather than in the two-thread configuration, the encoding, the instances or the seeds. Those re-runs also ended slightly better in 7 of the 11 cells and identical in the other 4, never worse, with the largest improvement close to 1\,\%, so the lost occupancy did carry a small quality cost. We keep the original cells in the analysis and treat the re-runs as a diagnostic rather than as data, since excluding runs on the basis of a host artifact would bias the very quantity the check exists to protect, and we report the episode as a reminder that effective width has to be measured per run rather than assumed from the requested one.

\subsection{Quality per thread level}
\label{sec:tier1-quality}

\begin{figure}[ht]
\centering
\includegraphics[width=\linewidth]{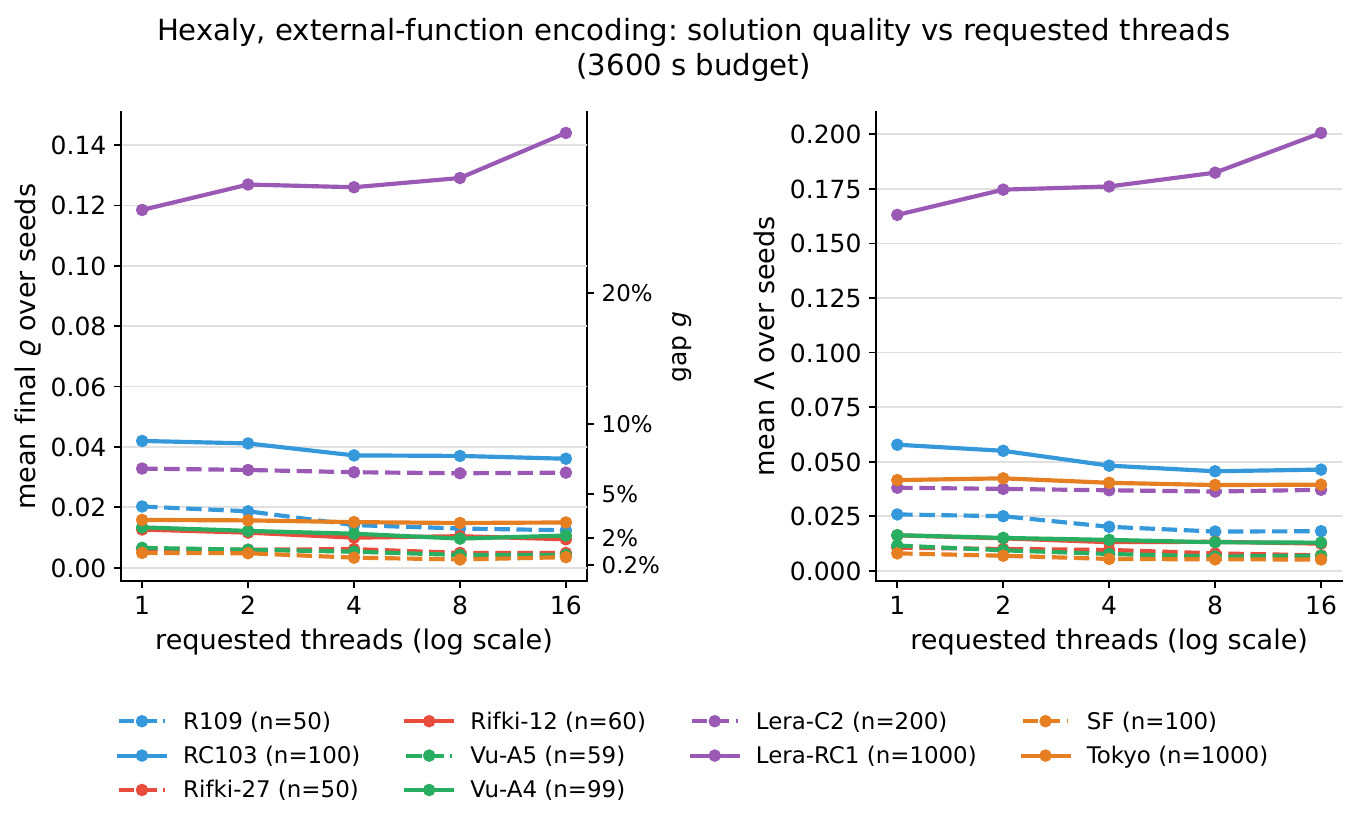}
\caption{External-function encoding under the C++ binding: solution quality per instance against the requested thread count, with ten seeds per point and the common one-hour budget. The left panel reports the mean final squeezed gap $\sgap$ and the right panel the mean anytime metric $\ascore$. The secondary axis of the left panel restates the final gap as the equivalent classic relative gap $g$. The anytime panel carries no such axis, for the reason given in \cref{sec:metrics}. Every point is a genuinely multi-threaded run at its nominal width. Nine of the ten instances end better at 16 threads than at one, while the tight-window instance of 1000 customers degrades and dominates the panel-equal mean.}
\label{fig:tier1-quality}
\end{figure}

Effective parallelism translates into very little pooled quality on this panel. Both metrics stay within a narrow band across the ladder and describe a shallow optimum at 8 threads (\cref{tab:pooled,fig:tier1-quality}). The anytime metric traces the same shape, with its best value 0.0368 at 8 threads. The widest setting therefore ends marginally ahead of the narrowest, by 0.6\,\% on the final gap and 0.7\,\% on the anytime metric, which is far too small a margin to read as a scaling effect. The mean time to the first checker-valid incumbent likewise moves little.

\begin{table}[ht]
\centering
\caption{Pooled results per thread level on the external-function encoding under the C++ binding (100 runs each: 10 instances $\times$ 10 seeds; seeds pooled per instance first, instances weighted equally). The dispersion column is the mean over instances of the sample standard deviation of final $\sgap$ across the ten seeds. The CPU/wall column tracks the requested width, so these are five genuinely distinct configurations. The pooled quality columns describe a shallow optimum at 8 threads, while the dispersion column falls until 8 threads and then rises.}
\label{tab:pooled}
\small
\begin{tabular}{rccccc}
\toprule
$\xi$ & mean final $\sgap$ & mean $\ascore$ & mean seed s.d. of $\sgap$ & mean first valid (s) & mean CPU/wall \\
\midrule
1 & 0.0273 & 0.0389 & 0.0042 & 9.5 & 1.00 \\
2 & 0.0275 & 0.0391 & 0.0036 & 9.5 & 1.85 \\
4 & 0.0260 & 0.0372 & 0.0033 & 9.1 & 3.84 \\
8 & 0.0257 & 0.0368 & 0.0029 & 8.8 & 7.77 \\
16 & 0.0271 & 0.0386 & 0.0032 & 8.8 & 15.75 \\
\bottomrule
\end{tabular}
\end{table}

The pooled band hides a split that is now cleanly attributable to a single instance. Nine of the ten instances end better at 16 threads than at one, and their anytime metric improves as well (\cref{tab:tier1-gprime-matrix,tab:tier1-pip-matrix}). The exception is the tight-window Lera-RC1 instance, which degrades from 0.1185 at one thread to 0.1440 at 16, an increase of 21.5\,\%, and it degrades on the anytime metric by a comparable 23.0\,\%. Since instances carry equal weight, that single degradation cancels the gains of the other nine. Indeed, pooling the same ladder over the nine instances other than Lera-RC1 gives a monotone descent of 17.5\,\% across the ladder. The pooled rows therefore do not show that additional cores are useless for this encoding, but that their benefit does not reach the hardest instance of the panel, where a wider search appears to spend its extra capacity without converting it into a better incumbent.

\begin{table}[ht]
\centering
\caption{Mean final $\sgap$ by instance and thread level on the external-function encoding under the C++ binding (mean over the ten seeds). The last column is the change from the $\xi = 1$ to the $\xi = 16$ thread-count entry, as a signed percentage of the $\xi = 1$ entry. Negative values mean the wider configuration ends better.}
\label{tab:tier1-gprime-matrix}
\small
\begin{tabular}{lrcccccc}
\toprule
Instance & $n$ & $\xi{=}1$ & $\xi{=}2$ & $\xi{=}4$ & $\xi{=}8$ & $\xi{=}16$ & $\Delta_{16-1}$ (\%) \\
\midrule
R109 & 50 & 0.0203 & 0.0187 & 0.0141 & 0.0130 & 0.0124 & $-38.9\%$ \\
RC103 & 100 & 0.0420 & 0.0412 & 0.0372 & 0.0370 & 0.0361 & $-14.0\%$ \\
Rifki-27 & 50 & 0.0062 & 0.0059 & 0.0062 & 0.0049 & 0.0049 & $-21.8\%$ \\
Rifki-12 & 60 & 0.0126 & 0.0116 & 0.0099 & 0.0104 & 0.0094 & $-25.6\%$ \\
Vu-A5 & 59 & 0.0066 & 0.0059 & 0.0053 & 0.0042 & 0.0042 & $-36.3\%$ \\
Vu-A4 & 99 & 0.0134 & 0.0122 & 0.0113 & 0.0096 & 0.0106 & $-20.3\%$ \\
Lera-C2 & 200 & 0.0328 & 0.0324 & 0.0316 & 0.0313 & 0.0315 & $-4.1\%$ \\
Lera-RC1 & 1000 & 0.1185 & 0.1269 & 0.1260 & 0.1291 & 0.1440 & $+21.5\%$ \\
SF & 100 & 0.0049 & 0.0048 & 0.0033 & 0.0027 & 0.0035 & $-28.7\%$ \\
Tokyo & 1000 & 0.0158 & 0.0157 & 0.0151 & 0.0148 & 0.0150 & $-5.3\%$ \\
\bottomrule
\end{tabular}
\end{table}

\begin{table}[ht]
\centering
\caption{Mean $\ascore$ by instance and thread level on the external-function encoding under the C++ binding (mean over the ten seeds). The anytime metric agrees with the final gap on the direction of every instance of the panel.}
\label{tab:tier1-pip-matrix}
\small
\begin{tabular}{lrcccccc}
\toprule
Instance & $n$ & $\xi{=}1$ & $\xi{=}2$ & $\xi{=}4$ & $\xi{=}8$ & $\xi{=}16$ & $\Delta_{16-1}$ (\%) \\
\midrule
R109 & 50 & 0.0258 & 0.0250 & 0.0202 & 0.0180 & 0.0181 & $-29.8\%$ \\
RC103 & 100 & 0.0578 & 0.0550 & 0.0482 & 0.0456 & 0.0464 & $-19.7\%$ \\
Rifki-27 & 50 & 0.0106 & 0.0100 & 0.0096 & 0.0080 & 0.0071 & $-33.4\%$ \\
Rifki-12 & 60 & 0.0162 & 0.0149 & 0.0131 & 0.0133 & 0.0124 & $-23.4\%$ \\
Vu-A5 & 59 & 0.0116 & 0.0094 & 0.0077 & 0.0066 & 0.0069 & $-40.5\%$ \\
Vu-A4 & 99 & 0.0163 & 0.0151 & 0.0142 & 0.0131 & 0.0129 & $-20.9\%$ \\
Lera-C2 & 200 & 0.0381 & 0.0375 & 0.0369 & 0.0363 & 0.0371 & $-2.5\%$ \\
Lera-RC1 & 1000 & 0.1631 & 0.1746 & 0.1761 & 0.1824 & 0.2006 & $+23.0\%$ \\
SF & 100 & 0.0080 & 0.0069 & 0.0054 & 0.0053 & 0.0052 & $-35.0\%$ \\
Tokyo & 1000 & 0.0416 & 0.0424 & 0.0403 & 0.0393 & 0.0395 & $-5.0\%$ \\
\bottomrule
\end{tabular}
\end{table}

\subsection{Dispersion across seeds}
\label{sec:tier1-dispersion}

\begin{figure}[ht]
\centering
\includegraphics[width=\linewidth]{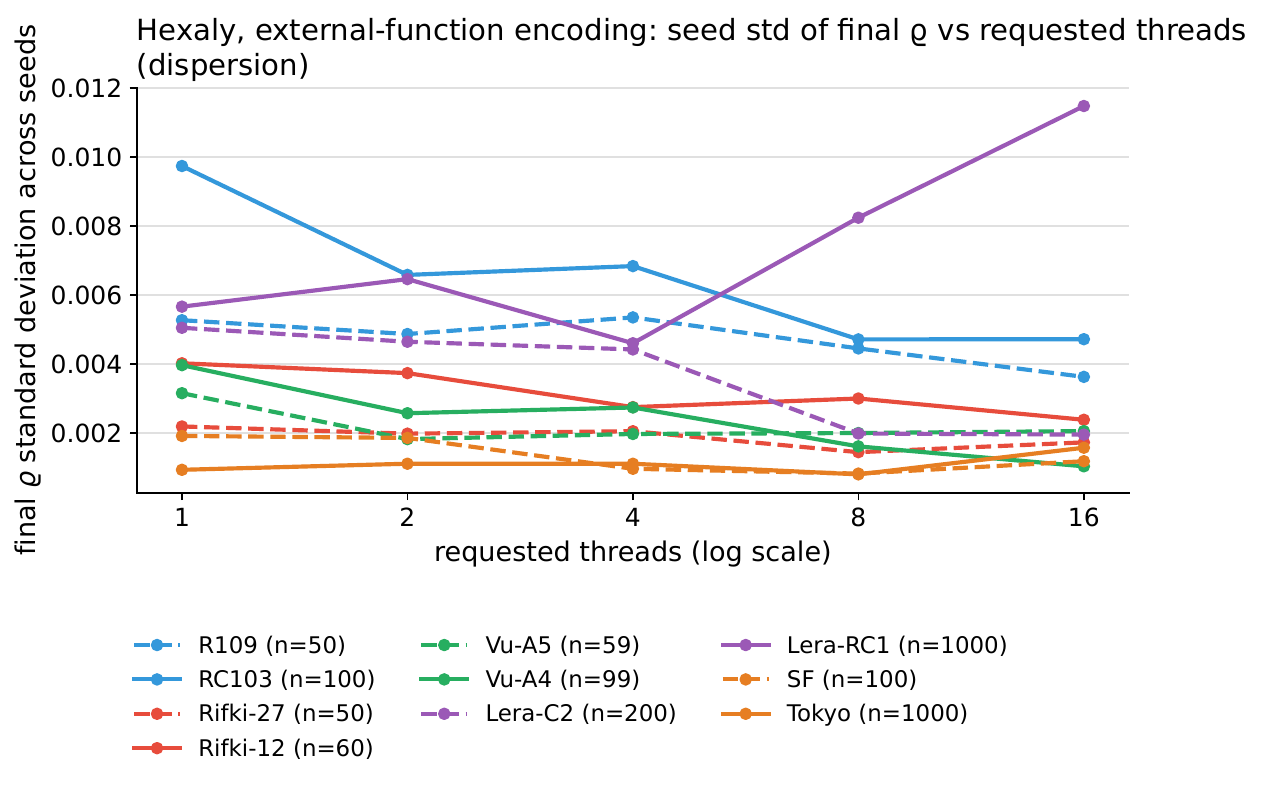}
\caption{External-function encoding under the C++ binding: sample standard deviation of final $\sgap$ across the ten seeds, per instance, against the requested thread count. The quantity lives on the squeezed scale and has no unique counterpart on the classic gap scale, so this figure carries a single panel and no gap axis. Every instance records a valid incumbent in every cell, so no cell is dispersed by a missing solution.}
\label{fig:tier1-std}
\end{figure}

Stability responds to width even where mean quality does not. The pooled seed standard deviation of the final gap falls by 30.6\,\% between 1 and 8 threads, then rises again at 16, where it remains 24.2\,\% below its one-thread value (\cref{tab:pooled,tab:tier1-stds,fig:tier1-std}). Under the Python binding the same quantity was flat at every width, so this is the first measurement of the vendor's stabilization account on this encoding. Eight of the ten instances are less dispersed at 16 threads than at one, and the reversal is carried by the tight-window instance, whose standard deviation at 16 threads is the largest of the ladder and 1.4 times its own value at 8 threads. Ten seeds support a standard deviation but not a claim about a population, and the panel is deliberately selected, so this remains a description of the panel rather than an estimated effect.

\begin{table}[ht]
\centering
\caption{Sample standard deviation across the ten seeds of final $\sgap$, by instance and thread level, on the external-function encoding under the C++ binding. The last column is the change from $\xi = 1$ to $\xi = 16$, as a signed percentage of the $\xi = 1$ entry. The percentages are large because the standard deviations themselves are small. The pooled mean of each column is the dispersion column of \cref{tab:pooled}. The widening at sixteen threads is carried by the Lera-RC1 row.}
\label{tab:tier1-stds}
\small
\begin{tabular}{lrcccccc}
\toprule
Instance & $n$ & $\xi{=}1$ & $\xi{=}2$ & $\xi{=}4$ & $\xi{=}8$ & $\xi{=}16$ & $\Delta_{16-1}$ (\%) \\
\midrule
R109 & 50 & 0.0053 & 0.0049 & 0.0054 & 0.0045 & 0.0036 & $-31.1\%$ \\
RC103 & 100 & 0.0097 & 0.0066 & 0.0068 & 0.0047 & 0.0047 & $-51.5\%$ \\
Rifki-27 & 50 & 0.0022 & 0.0020 & 0.0021 & 0.0015 & 0.0017 & $-21.0\%$ \\
Rifki-12 & 60 & 0.0040 & 0.0037 & 0.0028 & 0.0030 & 0.0024 & $-40.7\%$ \\
Vu-A5 & 59 & 0.0032 & 0.0018 & 0.0020 & 0.0020 & 0.0021 & $-34.8\%$ \\
Vu-A4 & 99 & 0.0040 & 0.0026 & 0.0027 & 0.0016 & 0.0010 & $-73.8\%$ \\
Lera-C2 & 200 & 0.0051 & 0.0046 & 0.0044 & 0.0020 & 0.0020 & $-61.2\%$ \\
Lera-RC1 & 1000 & 0.0057 & 0.0065 & 0.0046 & 0.0082 & 0.0115 & $+102.6\%$ \\
SF & 100 & 0.0019 & 0.0019 & 0.0010 & 0.0008 & 0.0012 & $-38.1\%$ \\
Tokyo & 1000 & 0.0009 & 0.0011 & 0.0011 & 0.0008 & 0.0016 & $+68.2\%$ \\
\bottomrule
\end{tabular}
\end{table}

\subsection{Memory and time to first solution}
\label{sec:tier1-cost}

Peak memory is essentially independent of the thread count on every instance, since the model data are shared and only per-worker state is replicated. The effect of width is visible only where the model itself is small, growing from about 31 to 58\,MB on the smallest instance, and it disappears at scale: the heavy urban instance holds 7.0\,GB at every width and the tight-window instance holds 0.6\,GB. Both are markedly cheaper than under the Python binding, which held 15.1\,GB and 1.8\,GB on the same two instances, so the port also relaxes the packing constraint that governs how many runs a node can host.

The first checker-valid incumbent appears within 10 seconds on average at every width, between 8.8 and 9.5 seconds across the ladder, against roughly 20 seconds under the Python binding. The pooled figure is dominated by the two largest instances, where data loading and model construction consume a visible part of the budget before the search begins.

\section{Choosing the slice count for the time-sliced encoding}
\label{sec:slices}

The time-sliced encoding introduces the slice count $\psi$, which controls the discretization of the travel-time horizon. We select it before the thread ladder because scaling results are meaningful only for a sufficiently faithful model. The single-threaded configuration study compares coarse, intermediate, and fine discretizations on six diverse instances under the same one-hour budget. Five instances subsequently appear in the thread ladder, so this is informed configuration selection rather than out-of-sample validation.

\subsection{The slice count affects checker feasibility}
\label{sec:slices-feasibility}

The main criterion for selecting the slice count is not the magnitude of the approximation error, but its effect on feasibility under the original travel-time functions. With a coarse discretization, a run may fail to record any incumbent accepted by the checker. \Cref{tab:slices-gprime} reports the mean final $\sgap$ for every instance and slice count, where $\sgap = 1$ denotes the absence of a checker-valid recorded incumbent.

\begin{table}[ht]
\centering
\caption{Time-sliced encoding: mean final $\sgap$ by instance and slice count (mean over ten seeds, 3600\,s budget, single-threaded). A value of 1 means no checker-valid incumbent was recorded. A parenthesis reports how many of the ten seeds recorded no checker-valid incumbent, and is shown only where that count is nonzero: every other cell had ten valid seeds. Those seeds enter the mean at $\sgap = 1$, so an annotated mean measures feasibility failure rather than solution quality. The last column is the change from $\psi = 5$ to $\psi = 96$, as a signed percentage of the $\psi = 5$ entry.}
\label{tab:slices-gprime}
\small
\begin{tabular}{lrr@{}lr@{}lrr}
\toprule
Instance & $n$ & \multicolumn{2}{c}{$\psi = 5$} & \multicolumn{2}{c}{$\psi = 24$} & $\psi = 96$ & $\Delta_{96-5}$ (\%) \\
\midrule
RC103 & 100 & 0.2424 & \,(2/10) & 0.1324 & \,(1/10) & 0.0374 & $-84.6\%$ \\
RC207 & 100 & 0.0533 &  & 0.0515 &  & 0.0439 & $-17.7\%$ \\
Rifki-27 & 50 & 0.0100 &  & 0.0124 &  & 0.0095 & $-5.5\%$ \\
Vu-A4 & 99 & 0.0158 &  & 0.0098 &  & 0.0108 & $-31.7\%$ \\
Lera-RC1 & 1000 & 1.0000 & \,(10/10) & 0.9362 & \,(9/10) & 0.1273 & $-87.3\%$ \\
Tokyo & 1000 & 0.0165 &  & 0.0136 &  & 0.0144 & $-12.6\%$ \\
\bottomrule
\end{tabular}
\end{table}

\begin{figure}[ht]
\centering
\includegraphics[width=\linewidth]{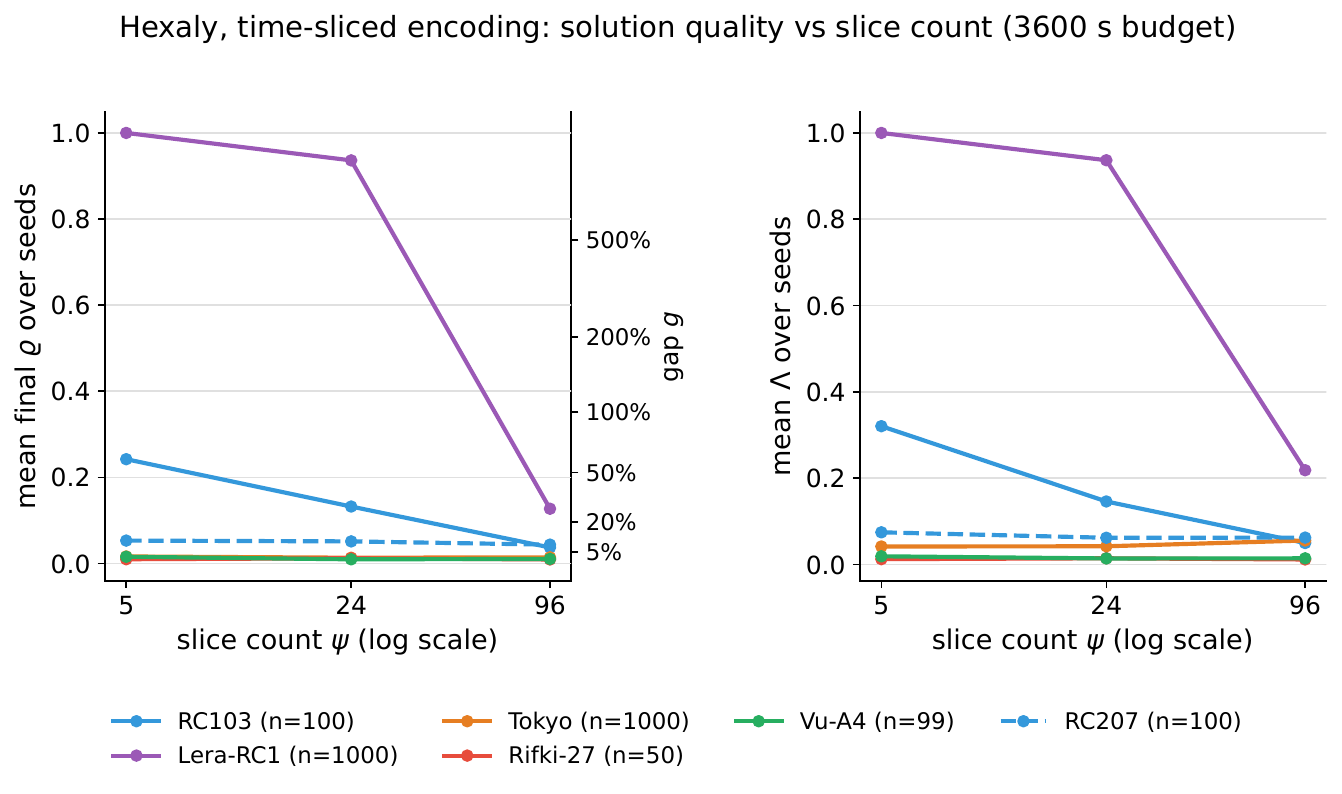}
\caption{Time-sliced encoding: solution quality per instance against the slice count. The left panel reports the mean final squeezed gap $\sgap$ and the right panel the mean anytime metric $\ascore$. The secondary axis of the left panel restates the final gap as the equivalent classic relative gap $g$. The anytime panel carries no such axis, for the reason given in \cref{sec:metrics}. Both metrics use the checker-valid recorded trajectory on the original travel-time functions. On the tight-window Lera-RC1 instance, $\sgap = 1$ at the coarsest setting denotes the absence of any checker-valid recorded incumbent rather than a large finite gap. The finest setting provides the most reliable feasibility and the strongest overall quality across both metrics.}
\label{fig:slices-quality}
\end{figure}

The tight-window Lera-RC1 instance provides the clearest example: coarse discretizations often produce no checker-valid recorded incumbent, whereas the finest setting succeeds for every seed (\cref{tab:slices-gprime}). Similar failures also occur on smaller instances. The discretized model can admit schedules that are infeasible under the original travel-time functions, which is why all solutions are evaluated through the reference checker. At $\psi = 96$, every run records a checker-valid incumbent and the checker rejects fewer trajectory points than under the coarser discretizations. Occasional invalid terminal solutions do not affect the final-gap metric, which uses the best valid recorded incumbent. \Cref{tab:slices-pooled,fig:slices-quality} summarize the results.

\begin{table}[ht]
\centering
\caption{Time-sliced encoding pooled per slice count (6 instances $\times$ 10 seeds; seeds pooled per instance first, instances weighted equally). Rejected incumbents are trajectory points the checker refused, per run, and are the fidelity readout of the discretization. The dispersion column carries, in parentheses, the number of runs of that row with no checker-valid incumbent, which enter the mean at $\sgap = 1$. Where that count is nonzero the dispersion measures how often the discretization fails rather than how variable its solutions are. Seed-to-seed stability of solution quality is therefore read from the thread ladders rather than from this table.}
\label{tab:slices-pooled}
\small
\begin{tabular}{rccccc}
\toprule
$\psi$ & mean final $\sgap$ & mean $\ascore$ & mean seed s.d. of $\sgap$ & mean first valid (s) & mean rejected incumbents \\
\midrule
5 & 0.2230 & 0.2448 & 0.0692 \,(12/60) & 73.7 & 361.4 \\
24 & 0.1927 & 0.2027 & 0.0887 \,(10/60) & 19.4 & 344.9 \\
96 & 0.0405 & 0.0688 & 0.0079 \,(0/60) & 56.4 & 186.8 \\
\bottomrule
\end{tabular}
\end{table}

Two readings of \cref{tab:slices-pooled} have to be kept apart. The pooled columns at $\psi = 5$ and $\psi = 24$ are dominated by the washout convention, since 12 and 10 of their 60 runs respectively contribute the worst attainable value $\sgap = 1$, thus their dispersion figures measure the frequency of that failure rather than the variability of a solution. Note that the two coarse settings therefore admit no reliable ordering against each other, and that the only comparison surviving the bookkeeping is the one that matters, namely that $\psi = 96$ is far ahead of both. Restricting the pool to the four instances where feasibility is not at stake, the fine discretization also wins by a margin an order of magnitude smaller than the pooled table suggests.

\subsection{Memory cost of \texorpdfstring{$\psi = 96$}{psi = 96}}
\label{sec:slices-cost}

The discretized table has size $\left| \mathcal{V} \right| \times \left| \mathcal{V} \right| \times \psi$, so memory grows linearly with the slice count and quadratically with the number of vertices. It is negligible on the small instances but reaches double-digit gigabytes on the largest one (\cref{tab:slices-rss}). Under the Python implementation the same runs peaked far higher, up to 25.0\,GB on the heaviest instance at $\psi = 96$, most of it interpreter-side table construction rather than the engine's own copy; the peaks reported here are what the engine itself holds. The retained configuration fits the available nodes with a wide margin, and its memory use changes little with thread count. On the largest instances, model construction consumes a visible but modest part of the wall-clock budget, slightly over two minutes on the heaviest instance, where the Python implementation spent more than four.

\begin{table}[!htbp]
\centering
\caption{Time-sliced encoding: peak resident set size (GNU~time maximum RSS, max over the ten seeds) by instance and slice count.}
\label{tab:slices-rss}
\small
\begin{tabular}{lrrrr}
\toprule
Instance & $n$ & $\psi = 5$ & $\psi = 24$ & $\psi = 96$ \\
\midrule
RC103 & 100 & 43\,MB & 47\,MB & 64\,MB \\
RC207 & 100 & 43\,MB & 47\,MB & 62\,MB \\
Rifki-27 & 50 & 42\,MB & 44\,MB & 50\,MB \\
Vu-A4 & 99 & 101\,MB & 100\,MB & 102\,MB \\
Lera-RC1 & 1000 & 0.6\,GB & 1.1\,GB & 2.8\,GB \\
Tokyo & 1000 & 7.0\,GB & 7.0\,GB & 12.7\,GB \\
\bottomrule
\end{tabular}
\end{table}

In summary, the lower memory cost of $\psi = 5$ and $\psi = 24$ does not compensate for their occasional failure to produce checker-valid recorded incumbents. The $\psi = 96$ configuration records at least one valid incumbent for every instance and seed, gives the best or a near-best final quality throughout the study, and remains within the available memory budget. We consequently retain $\psi = 96$ for the thread ladder of \cref{sec:tier2}.

\section{The thread ladder on the time-sliced encoding}
\label{sec:tier2}

\subsection{Setup}
\label{sec:tier2-setup}

This second ladder measures effective thread scaling. It reuses the instances, seeds, thread levels, wall-clock budget, and hardware class of the external-function ladder, enabling paired comparisons throughout, and fixes the slice count at $\psi = 96$ according to \cref{sec:slices}. Runs follow the execution and scoring protocol of \cref{sec:grid,sec:metrics}: each process is pinned to exclusive physical cores, monitored independently, and scored by the reference checker. Nodes are grouped by thread width and packed within their memory capacity. This prevents core overlap but leaves width partly confounded with co-tenancy, since narrow configurations place more processes on a node. Every run records a checker-valid incumbent. A small number of terminal solutions on the hardest instance are rejected, but no run is excluded because the metrics use the best checker-valid recorded incumbent.

\subsection{Quality per thread level}
\label{sec:tier2-quality}

\begin{figure}[ht]
\centering
\includegraphics[width=\linewidth]{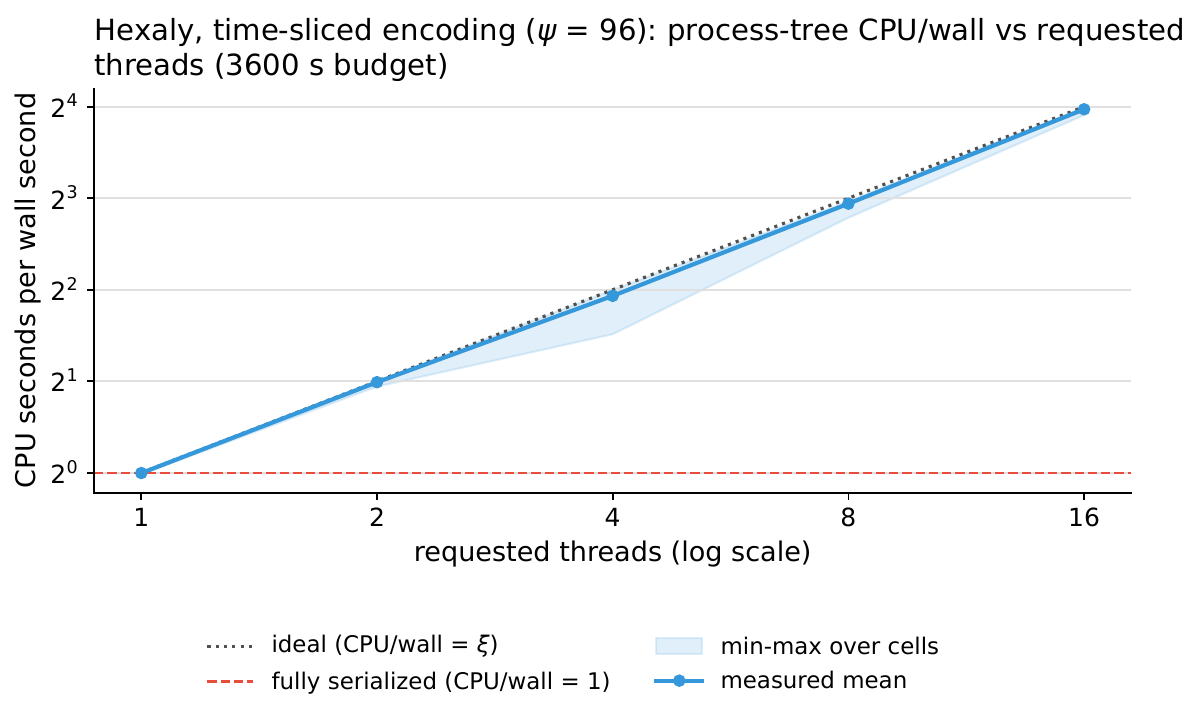}
\caption{Time-sliced encoding: measured process-tree CPU per wall-clock second against the requested thread count (GNU~time, whole process tree, 100 runs per level, mean and min--max band). The diagonal denotes full occupancy of the allocated cores. This figure shows similar occupancy to \cref{fig:tier1-cpu}.}
\label{fig:tier2-cpu}
\end{figure}

The time-sliced encoding keeps the allocated cores close to fully occupied at every width (\cref{fig:tier2-cpu}), as the external-function encoding does under its C++ binding. A modest dip at the intermediate width persists, with a low tail on a few of the smaller instances, and it reproduces at the same width in the superseded Python-binding runs of the same cells, so it is a property of the engine at this width rather than of either binding. Neither preprocessing, CPU accounting, nor pinning explains it, and the available traces lack per-worker activity, so its solver-internal cause remains unidentified. Both encodings thus receive genuine parallelism, and what remains to be compared is what each does with it.

The additional cores yield substantial quality gains (\cref{tab:tier2-pooled}). From 1 to 16 threads, the pooled final gap falls by about 23\,\%, and the anytime metric improves in 99 of the 100 paired runs. The final gap improves in 93 of them and on every per-instance mean. Quality continues to progress at the widest measured level, so no saturation is observed within this ladder. This is a wall-time comparison and the widest setting uses 16 times the nominal process width.

\begin{table}[ht]
\centering
\caption{Pooled results per thread level on the time-sliced encoding at $\psi = 96$ (100 runs each: 10 instances $\times$ 10 seeds; seeds pooled per instance first, instances weighted equally). The dispersion column is the mean over instances of the sample standard deviation of final $\sgap$ across the ten seeds. The CPU/wall ratio follows the requested thread count here as it does for the external-function encoding in \cref{tab:pooled}. Unlike that ladder, both quality columns improve monotonically with width, and the dispersion column ends at nearly half its narrow-width value.}
\label{tab:tier2-pooled}
\small
\begin{tabular}{rccccc}
\toprule
$\xi$ & mean final $\sgap$ & mean $\ascore$ & mean seed s.d. of $\sgap$ & mean first valid (s) & mean CPU/wall \\
\midrule
1 & 0.0275 & 0.0462 & 0.0056 & 38.9 & 1.00 \\
2 & 0.0261 & 0.0441 & 0.0038 & 35.3 & 1.99 \\
4 & 0.0237 & 0.0413 & 0.0032 & 26.4 & 3.82 \\
8 & 0.0230 & 0.0373 & 0.0034 & 24.9 & 7.68 \\
16 & 0.0211 & 0.0345 & 0.0031 & 22.3 & 15.69 \\
\bottomrule
\end{tabular}
\end{table}

\begin{figure}[ht]
\centering
\includegraphics[width=\linewidth]{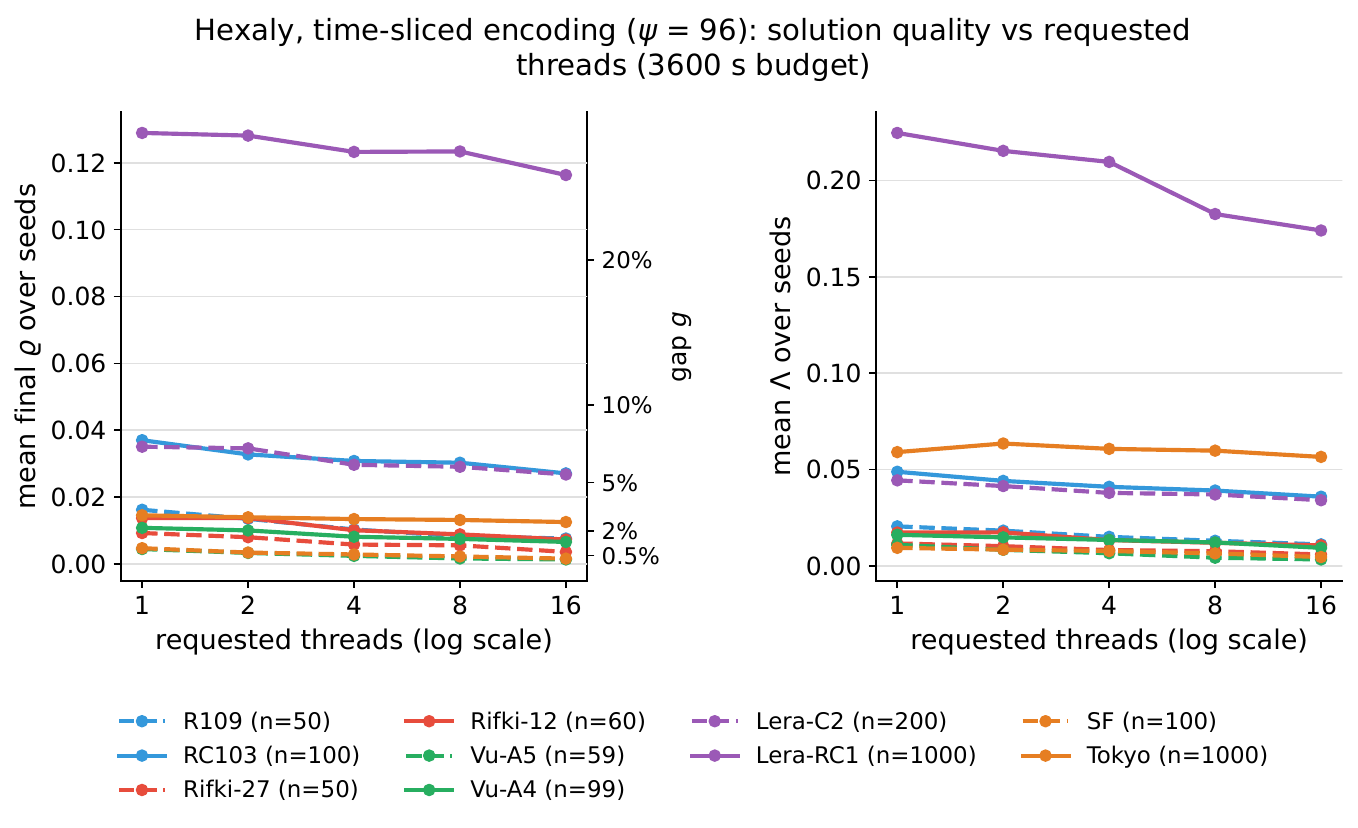}
\caption{Time-sliced encoding at $\psi = 96$: solution quality per instance against the thread count, with ten seeds per point and the common one-hour budget. The left panel reports the mean final squeezed gap $\sgap$ and the right panel the mean anytime metric $\ascore$. The secondary axis of the left panel restates the final gap as the equivalent classic relative gap $g$. The anytime panel carries no such axis, for the reason given in \cref{sec:metrics}. Every instance ends better at 16 threads than at one, and the anytime metric shows the same broad benefit from additional cores. Intermediate widths need not improve monotonically, but both panels support a substantial overall gain across the tested panel.}
\label{fig:tier2-quality}
\end{figure}

The gain is broad rather than confined to the largest instances: every instance ends better at 16 threads than at one (\cref{tab:tier2-gprime-matrix,tab:tier2-pip-matrix,fig:tier2-quality}). Intermediate levels are not always monotone, particularly on the hardest instance, but the overall direction is consistent across sizes.

\begin{table}[ht]
\centering
\caption{Mean final $\sgap$ by instance and thread level on the time-sliced encoding at $\psi = 96$ (mean over the ten seeds). The last column is the change from $\xi = 1$ to $\xi = 16$, as a signed percentage of the $\xi = 1$ entry.}
\label{tab:tier2-gprime-matrix}
\small
\begin{tabular}{lrcccccc}
\toprule
Instance & $n$ & $\xi{=}1$ & $\xi{=}2$ & $\xi{=}4$ & $\xi{=}8$ & $\xi{=}16$ & $\Delta_{16-1}$ (\%) \\
\midrule
R109 & 50 & 0.0162 & 0.0135 & 0.0103 & 0.0083 & 0.0076 & $-53.1\%$ \\
RC103 & 100 & 0.0370 & 0.0328 & 0.0308 & 0.0303 & 0.0271 & $-26.8\%$ \\
Rifki-27 & 50 & 0.0092 & 0.0080 & 0.0058 & 0.0056 & 0.0036 & $-61.1\%$ \\
Rifki-12 & 60 & 0.0137 & 0.0137 & 0.0101 & 0.0088 & 0.0073 & $-46.6\%$ \\
Vu-A5 & 59 & 0.0045 & 0.0033 & 0.0024 & 0.0016 & 0.0013 & $-70.0\%$ \\
Vu-A4 & 99 & 0.0108 & 0.0100 & 0.0081 & 0.0075 & 0.0066 & $-39.1\%$ \\
Lera-C2 & 200 & 0.0351 & 0.0346 & 0.0297 & 0.0291 & 0.0267 & $-23.7\%$ \\
Lera-RC1 & 1000 & 0.1289 & 0.1281 & 0.1232 & 0.1234 & 0.1163 & $-9.8\%$ \\
SF & 100 & 0.0047 & 0.0034 & 0.0028 & 0.0022 & 0.0016 & $-66.3\%$ \\
Tokyo & 1000 & 0.0146 & 0.0140 & 0.0134 & 0.0132 & 0.0125 & $-13.9\%$ \\
\bottomrule
\end{tabular}
\end{table}

\begin{table}[!htbp]
\centering
\caption{Mean $\ascore$ by instance and thread level on the time-sliced encoding at $\psi = 96$ (mean over the ten seeds). The anytime metric confirms the direction of every row of \cref{tab:tier2-gprime-matrix}.}
\label{tab:tier2-pip-matrix}
\small
\begin{tabular}{lrcccccc}
\toprule
Instance & $n$ & $\xi{=}1$ & $\xi{=}2$ & $\xi{=}4$ & $\xi{=}8$ & $\xi{=}16$ & $\Delta_{16-1}$ (\%) \\
\midrule
R109 & 50 & 0.0204 & 0.0182 & 0.0150 & 0.0129 & 0.0111 & $-45.8\%$ \\
RC103 & 100 & 0.0488 & 0.0440 & 0.0410 & 0.0390 & 0.0358 & $-26.5\%$ \\
Rifki-27 & 50 & 0.0116 & 0.0101 & 0.0081 & 0.0075 & 0.0058 & $-49.7\%$ \\
Rifki-12 & 60 & 0.0172 & 0.0173 & 0.0133 & 0.0119 & 0.0105 & $-38.8\%$ \\
Vu-A5 & 59 & 0.0108 & 0.0082 & 0.0065 & 0.0041 & 0.0032 & $-70.2\%$ \\
Vu-A4 & 99 & 0.0160 & 0.0147 & 0.0134 & 0.0121 & 0.0093 & $-42.2\%$ \\
Lera-C2 & 200 & 0.0443 & 0.0413 & 0.0378 & 0.0370 & 0.0340 & $-23.3\%$ \\
Lera-RC1 & 1000 & 0.2247 & 0.2154 & 0.2096 & 0.1826 & 0.1740 & $-22.5\%$ \\
SF & 100 & 0.0093 & 0.0083 & 0.0076 & 0.0064 & 0.0045 & $-51.6\%$ \\
Tokyo & 1000 & 0.0590 & 0.0634 & 0.0607 & 0.0597 & 0.0565 & $-4.2\%$ \\
\bottomrule
\end{tabular}
\end{table}

The gain is not explained by a larger reported search volume: multi-threaded runs report slightly fewer iterations than the single-threaded baseline. We use this counter only as a diagnostic because the documentation does not say whether it aggregates worker activity or counts synchronized search rounds. The result therefore rules out a simple increase in the reported iteration count, but it cannot establish a diversification mechanism. It is this monotone response to width, rather than any advantage at a fixed width, that separates the two encodings, as developed in \cref{sec:headtohead}.

\subsection{Dispersion across seeds}
\label{sec:tier2-dispersion}

\begin{figure}[!htbp]
\centering
\includegraphics[width=\linewidth]{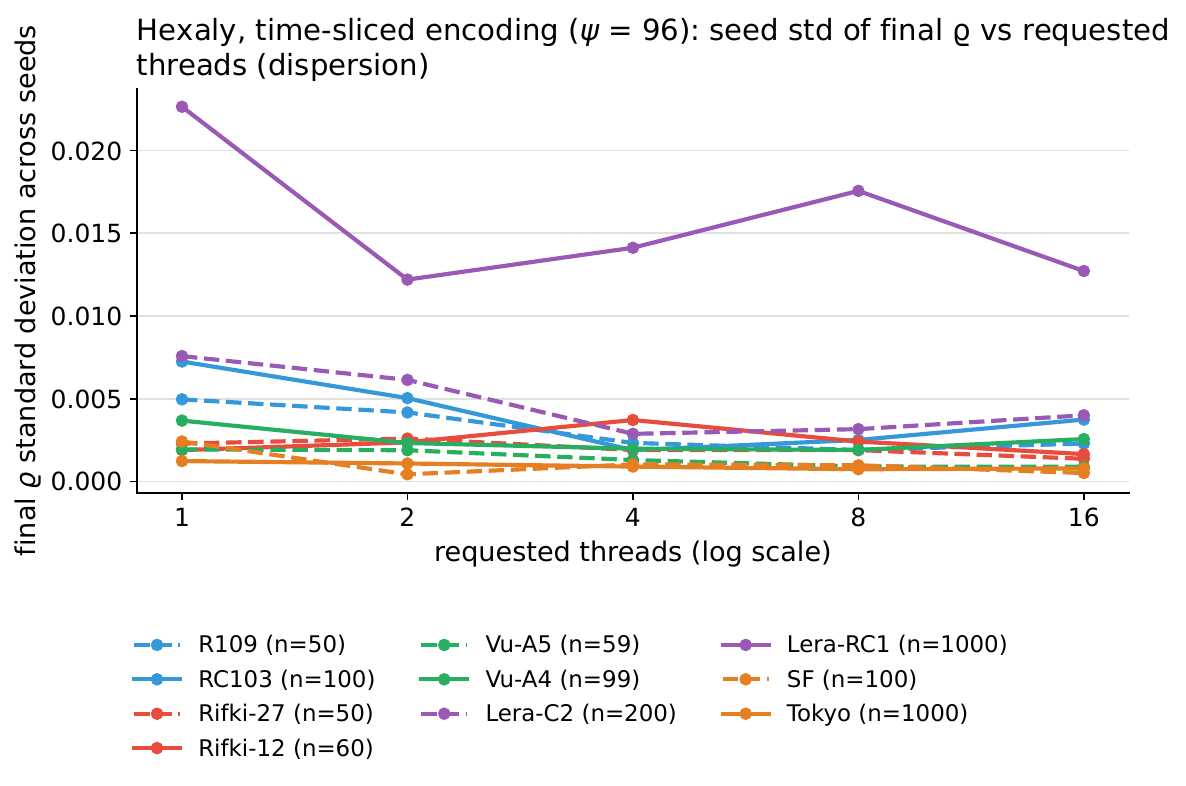}
\caption{Time-sliced encoding at $\psi = 96$: sample standard deviation of final $\sgap$ across the ten seeds, per instance, against the thread count. The pooled standard deviation at 16 threads is close to half its one-thread value, and every instance of the panel moves in that direction. As in \cref{fig:tier1-std}, the quantity lives on the squeezed scale alone, so the figure carries a single panel.}
\label{fig:tier2-std}
\end{figure}

The seed standard deviations of this encoding provide a second descriptive measure of stability. The pooled standard deviation drops sharply over the first half of the ladder and ends with a reduction of about 45\,\%, with all ten instances being less dispersed at the wide end than at the narrow one (\cref{tab:tier2-stds,fig:tier2-std}).

\begin{table}[ht]
\centering
\caption{Sample standard deviation across the ten seeds of final $\sgap$, by instance and thread level, on the time-sliced encoding at $\psi = 96$. The last column is the change from $\xi = 1$ to $\xi = 16$, as a signed percentage of the $\xi = 1$ entry. Negative values mean the wider configuration is more stable across seeds, which holds for every instance here.}
\label{tab:tier2-stds}
\small
\begin{tabular}{lrcccccc}
\toprule
Instance & $n$ & $\xi{=}1$ & $\xi{=}2$ & $\xi{=}4$ & $\xi{=}8$ & $\xi{=}16$ & $\Delta_{16-1}$ (\%) \\
\midrule
R109 & 50 & 0.0050 & 0.0042 & 0.0023 & 0.0019 & 0.0023 & $-53.8\%$ \\
RC103 & 100 & 0.0072 & 0.0050 & 0.0019 & 0.0025 & 0.0037 & $-48.4\%$ \\
Rifki-27 & 50 & 0.0023 & 0.0026 & 0.0019 & 0.0019 & 0.0014 & $-40.0\%$ \\
Rifki-12 & 60 & 0.0019 & 0.0024 & 0.0037 & 0.0024 & 0.0017 & $-12.8\%$ \\
Vu-A5 & 59 & 0.0019 & 0.0019 & 0.0013 & 0.0009 & 0.0009 & $-54.1\%$ \\
Vu-A4 & 99 & 0.0037 & 0.0023 & 0.0020 & 0.0019 & 0.0026 & $-30.3\%$ \\
Lera-C2 & 200 & 0.0076 & 0.0061 & 0.0029 & 0.0032 & 0.0040 & $-47.2\%$ \\
Lera-RC1 & 1000 & 0.0227 & 0.0122 & 0.0141 & 0.0176 & 0.0127 & $-43.8\%$ \\
SF & 100 & 0.0024 & 0.0004 & 0.0010 & 0.0010 & 0.0005 & $-78.8\%$ \\
Tokyo & 1000 & 0.0012 & 0.0011 & 0.0009 & 0.0007 & 0.0008 & $-35.9\%$ \\
\bottomrule
\end{tabular}
\end{table}

Worst-seed outcomes improve together with the dispersion. On the most difficult instance, additional threads mitigate an unfavorable seed more strongly than they improve the mean, which is valuable when selecting a practical configuration. The endpoint standard deviations are descriptively consistent with the stabilization mechanism described by Hexaly on this panel. Ten seeds give a usable standard deviation per cell but the instances are deliberately selected, so nothing here justifies generalization beyond the panel. The external-function encoding narrows its dispersion as well, by about 30\,\% up to 8 threads before giving part of it back at 16 (\cref{sec:tier1-dispersion}), so stabilization appears in both encodings and is more robust than the quality response, which appears in only one.

\subsection{Memory and time to first solution}
\label{sec:tier2-cost}

Peak memory remains nearly independent of thread count across the panel. A wider run therefore consumes additional cores but little additional memory, and the budget established in \cref{sec:slices-cost} remains applicable. The first checker-valid incumbent appears earlier at every step of the ladder on the pooled mean, from 38.9\,s at one thread to 22.3\,s at 16, driven mainly by the large instances. The level itself sits well below the corresponding Python-implementation column, whose one-thread mean was 64.1\,s, because the C++ implementation builds the sliced model faster and starts searching sooner. As such, that shift is a property of the implementations, not of width. On the difficult Lera instance the earlier valid solution accompanies genuine additional CPU use and is consistent with a broader parallel search.

\subsection{Comparison of the two encodings at equal budget}
\label{sec:headtohead}

The two ladders share their instances, seeds, thread levels, wall-clock budget, hardware class and checker, so they can be compared cell by cell. At one thread the encodings are close, with the exact one being ahead by only 0.6\,\% on the pooled mean, a margin small enough that its composition matters more than its size (\cref{tab:headtohead,tab:headtohead-pooled}). Indeed, the exact encoding wins only 4 of the 10 instances at that width, and its entire pooled advantage comes from the tight-window instance where the discretization struggles. Removing that one instance reverses the sign and leaves the sliced encoding 5.7\,\% ahead over the remaining nine. The discretization therefore buys no clear quality per unit of time on this panel once the exact evaluation is not paying an interface tax, which contradicts the uniform advantage reported in the first version of this report which previously measured the cost of the Python bridge rather than a property of the encodings.

The measurement is binding-symmetric: both arms reach the solver through its C++ interface, and the per-run parity check of \cref{sec:model} shows they solved identical sliced tables cell for cell, so what is compared is the encodings rather than the interfaces. An earlier form of this comparison still charged the sliced arm its Python-side model construction inside the budget, which delayed its first valid incumbent on the largest instances (45 to 64 seconds against roughly 9 for the exact arm) and put the exact edge at one thread at 1.5\,\%. Removing that charge improves the sliced arm's anytime metric by 10 to 17\,\% relative at every width while leaving its final gaps essentially unchanged, and narrows the one-thread edge to its present 0.6\,\%: the residual is the encoding difference, not an interface artifact.

Width reverses the ordering. The sliced encoding improves at every step while the exact one does not, so the pooled mean crosses between 1 and 2 threads and the margin then widens at every further step, reaching 22.4\,\% at 16 threads. Across the 50 pairs of instance and thread level, the exact encoding wins 10, and at the widest setting it wins none. The practical advantage of the discretization is thus thread scaling, not a better exchange of fidelity for throughput.

\begin{table}[ht]
\centering
\caption{Mean final $\sgap$ of the two encodings at the narrowest and widest thread level, 10 seeds, one-hour budget, both scored by the checker on the original travel-time functions. Exact denotes the external-function encoding under the C++ binding and sliced denotes the time-sliced encoding at $\psi = 96$. Each $\Delta$ is the sliced entry read against the exact one, as a signed percentage of the exact entry. A negative $\Delta$ means the sliced encoding ends ahead. The last row pools instances with equal weight, as everywhere in this report. Intermediate thread levels are omitted here and pooled in \cref{tab:headtohead-pooled}.}
\label{tab:headtohead}
\small
\begin{tabular}{lrccrccr}
\toprule
& & \multicolumn{3}{c}{$\xi = 1$} & \multicolumn{3}{c}{$\xi = 16$} \\
\cmidrule(lr){3-5}\cmidrule(lr){6-8}
Instance & $n$ & exact & sliced & $\Delta$ & exact & sliced & $\Delta$ \\
\midrule
R109 & 50 & 0.0203 & 0.0162 & $-20.1\%$ & 0.0124 & 0.0076 & $-38.6\%$ \\
RC103 & 100 & 0.0420 & 0.0370 & $-11.8\%$ & 0.0361 & 0.0271 & $-25.0\%$ \\
Rifki-27 & 50 & 0.0062 & 0.0092 & $+49.0\%$ & 0.0049 & 0.0036 & $-26.0\%$ \\
Rifki-12 & 60 & 0.0126 & 0.0137 & $+8.4\%$ & 0.0094 & 0.0073 & $-22.2\%$ \\
Vu-A5 & 59 & 0.0066 & 0.0045 & $-32.1\%$ & 0.0042 & 0.0013 & $-68.0\%$ \\
Vu-A4 & 99 & 0.0134 & 0.0108 & $-19.1\%$ & 0.0106 & 0.0066 & $-38.1\%$ \\
Lera-C2 & 200 & 0.0328 & 0.0351 & $+6.8\%$ & 0.0315 & 0.0267 & $-15.0\%$ \\
Lera-RC1 & 1000 & 0.1185 & 0.1289 & $+8.8\%$ & 0.1440 & 0.1163 & $-19.2\%$ \\
SF & 100 & 0.0049 & 0.0047 & $-3.8\%$ & 0.0035 & 0.0016 & $-54.4\%$ \\
Tokyo & 1000 & 0.0158 & 0.0146 & $-8.0\%$ & 0.0150 & 0.0125 & $-16.4\%$ \\
\midrule
pooled &  & 0.0273 & 0.0275 & $+0.6\%$ & 0.0271 & 0.0211 & $-22.4\%$ \\
\bottomrule
\end{tabular}
\end{table}

\begin{table}[ht]
\centering
\caption{Pooled mean final $\sgap$ of the two encodings at every thread level (10 instances $\times$ 10 seeds per cell, instances weighted equally). $\Delta$ is the sliced entry as a signed percentage of the exact one. The wins column counts the instances at that level where the exact encoding ends ahead. The two encodings cross between one and two threads. Only the sliced one improves with width.}
\label{tab:headtohead-pooled}
\small
\begin{tabular}{rcccc}
\toprule
$\xi$ & exact & sliced & $\Delta$ & exact wins \\
\midrule
1 & 0.0273 & 0.0275 & $+0.6\%$ & 4 of 10 \\
2 & 0.0275 & 0.0261 & $-5.0\%$ & 4 of 10 \\
4 & 0.0260 & 0.0237 & $-9.0\%$ & 1 of 10 \\
8 & 0.0257 & 0.0230 & $-10.5\%$ & 1 of 10 \\
16 & 0.0271 & 0.0211 & $-22.4\%$ & 0 of 10 \\
\bottomrule
\end{tabular}
\end{table}

Two further qualifications concern its scope. Both encodings are scored by the best checker-valid point of their recorded trajectory (\cref{sec:metrics}), a rule that matters more for the approximation, which produces many rejected incumbents, than for the exact model, whose native objective nearly coincides with the checker's. The comparison therefore concerns the encodings as deployed with independent validation, not their native objectives alone. And 10 instances with 10 seeds describe a panel rather than a population, so the near parity at one thread should be read as the absence of a visible advantage rather than as a proven equality. Under the single-core convention of our main solver comparison, the exact encoding is the default configuration we retain for Hexaly.

\section{Conclusion}
\label{sec:conclusion}

This report studied the thread scaling of Hexaly Optimizer 15.0 on two TDVRPTW encodings, as support experiments for the main solver comparison of \kayros. The external-function encoding evaluates continuous travel times exactly, and it uses the cores it is given once its functions are declared through a compiled C++ binding. Through the Python binding, which marks such callbacks as coming from a GIL language, it does not, as highlighted by the first version of this report. On this panel the exact encoding gains stability from additional cores, reducing its seed dispersion by about 30\,\% up to 8 threads, but almost no pooled quality, because nine of its ten instances improve while the hardest one degrades by enough to cancel them.

The time-sliced encoding approximates travel times natively and responds to width at every step. A separate configuration study led us to retain $\psi = 96$: coarser discretizations can fail checker validation, while this setting reliably produces valid incumbents and remains within the available memory budget. With this encoding, increasing the configuration from 1 to 16 threads lowers the pooled final gap by about 23\,\%, improves nearly every paired endpoint, and cuts the seed dispersion by nearly half. At a single thread the two encodings are 
on par, with the exact one ahead by only 0.6\,\% on this panel. The small pooled edge of exact evaluation rests entirely on the hardest instance, so the discretization earns its place through thread scaling rather than through a favorable exchange of fidelity for throughput. Stabilization is the one effect both encodings share, and it is the more robust of the two, since it survives an encoding whose mean quality does not improve. It is consistent with the diversification account given by Hexaly.

What the design supports is descriptive. Ten seeds give a usable dispersion per cell but no claim about a population, the 10 instances are a deliberately chosen panel rather than a sample, and width remains partly confounded with node co-tenancy, although the independent single-threaded replication of \cref{sec:tier2-setup} bounds that effect well below the observed gain. Why a wide search should lose ground on precisely the instance where feasible schedules are scarcest is not answered by our traces, which record no per-worker activity, and it is reported as an observation rather than as a mechanism. None of this revises the single-core convention of the main solver comparison. It quantifies a separate tradeoff: better and more stable solutions at equal wall time, in exchange for up to 16 times the nominal core allocation and little additional memory.

Two practical consequences follow for work with Hexaly. Its parameter API reports the requested thread count and not the effective one, so the width a Hexaly run actually obtains has to be established by independent CPU accounting instead of read back from the configuration. And because Hexaly declares Python external functions as callbacks from a GIL language and then restricts the effective search width accordingly, a Hexaly model whose external functions are written in Python should be scheduled, budgeted and reported as single-threaded, whereas the same model reaching the solver through its C++ interface should not. The Python interface also charges its model construction to the run's own budget: on our largest instance it consumed more than four minutes of the hour before the search began, where the C++ interface spends about two, and the anytime metric of the time-sliced arm moved by 10 to 17\,\% relative when that charge was removed. Both bindings solve the same model, but they do not leave it the same time to search. A claim about a solver's parallelism is only as good as the accounting behind it.

\section*{Data availability}

The solver runs behind every table and figure of this report and of its earlier versions are deposited as a public dataset \autocite{rascoussierHexalyThreadScalingData2026}. Each run is a JSON record carrying the configuration, the timing, the solver's native statistics and the complete incumbent trajectory, together with the independent GNU~time output for the same run. The deposit also contains the frozen reference snapshot of 2026/08/06 used to score every gap\footnote{Its entry set hashes to \texttt{5e9e4891\allowbreak 69d87cde\allowbreak 8287700501d88ebe\allowbreak a1b60afd\allowbreak 3731b49e\allowbreak 17827712e34b3563}.}, so the reported gaps can be recomputed exactly rather than approximately. The instances are not duplicated there: they are distributed as checksummed artifacts through the MAMUT-routing benchmark platform \autocite{pichon:hal-05629810v1}, and each record identifies its instance by a path relative to that tree. The model implementations, the run driver and the analysis code are not released at this stage.

\section*{Acknowledgments}

We thank Maxime Rougier and the Hexaly support team for providing the academic license, for their prompt and detailed answers, and for the remarks on multi-threading that motivated this study. We also thank Romain Billot, Christine Solnon and Lina Fahed for supervising the PhD in which this work was conducted. Special thanks to Guillaume Beslon and the BIOTIC team for their support. The experiments presented in this report were carried out using the Grid'5000 testbed, supported by a scientific interest group hosted by Inria and including CNRS, RENATER, several universities and other organizations (see \url{https://www.grid5000.fr}). This work is funded by the French National Research Agency (ANR) as part of the MAMUT project, ANR-22-CE22-0016, \enquote{Machine learning And Matheuristics algorithms for Urban Transportation}.

\section*{Declaration of AI use}

Generative AI tools were used throughout the preparation of this work, including frontier models such as GPT-5.6 Sol (OpenAI) and Claude Fable and Sonnet 5 (Anthropic). They contributed to the experiment orchestration, the analysis tooling and the drafting of this report. The author reviewed and revised all outputs from these tools and takes full responsibility for the final content of the work.

\printbibliography

\end{document}